\documentclass[11pt]{article}
\usepackage{amssymb}
\usepackage{amsmath}
\usepackage{latexsym}
\usepackage{amsfonts}

\newtheorem{thm}{Theorem}[section]

\newtheorem{cor}{Corollary}[section]
\newtheorem{Def}{Definition}[section]
\newtheorem{rem}{Remark}[section]
\newtheorem{exm}{Example}[section]

\begin{document}
\title{\bf Multivariate {\emph{q}}-P\'{o}lya and\\ inverse {\emph{q}}-P\'{o}lya distributions}
\author{Charalambos A. Charalambides\\
{\em Department of Mathematics, University of Athens}, \\
{\em Panepistemiopolis, GR-15784 Athens, Greece}}
\footnotetext{E-mail address: ccharal@math.uoa.gr}
\date{}
\maketitle

\begin{abstract}
An urn containing specified numbers of balls of distinct ordered colors is considered. A multiple $q$-P\'{o}lya urn model is introduced by assuming that the probability of $q$-drawing a ball of a specific color from the urn varies geometrically, with rate $q$, both with the number of drawings and the number of balls of the specific color, together with the total number of balls of the preceded colors, drawn in the previous $q$-drawings. Then, the joint distributions of the numbers of balls of distinct colors drawn (a) in a specific number of $q$-drawings and (b) until the occurrence of a specific number of balls of a certain color, are derived. These two distributions turned out to be $q$-analogues of the multivariate P\'{o}lya and inverse P\'{o}lya distributions, respectively. Also, the limiting distributions of the multivariate $q$-P\'{o}lya and inverse $q$-P\'{o}lya distributions, as the initial total number of balls in the urn tends to infinity, are shown to be $q$-multinomial and negative $q$-multinomial distributions, respectively.
\end{abstract}
{\em AMS(2000) subject classification}. Primary 60C05, Secondary 05A30.

\section{Introduction}\label{sec1}
A $q$-P\'{o}lya urn model is introduced in Charalambides (2012, 2016) and $q$-P\'{o}lya and inverse $q$-P\'{o}lya distributions are studied. Moreover, their limiting distributions, as the number of balls in the urn, or the number of drawn balls of one of the two colors, tends to infinity are shown to be $q$-binomial and negative $q$-binomial distributions. Kupershmidt (2000) introduced a $q$-hypergeometric distribution and a $q$-P\'{o}lya distribution (under the name $q$-contagious distribution) and represented the corresponding random variable as a sum of two-valued dependent random variables. Kemp (2005) starting from the Chu-Vandermonde sum as a probability generating function obtained two $q$-confluent hypergeometric distributions. She also, deduced these distributions as steady-state birth and death Markov chains.

The aim of this article is to introduced and studied in detail multivariate $q$-P\'{o}lya and inverse $q$-P\'{o}lya distributions and also, examined some of their limiting distributions. Section 2 is devoted to $q$-multinomial convolutions. Precisely, multivariate $q$-Vandermonde and inverse $q$-Vandermonde formulae are presented. Also, a closely connected multivariate $q$-Cauchy formula is deduced. In section 3, the $q$-P\'{o}lya urn model is extended to a multiple $q$-P\'{o}lya urn model by considering successive $q$-drawings of balls from an urn containing specified numbers of balls of different colors and assuming that the probability of randomly $q$-drawing a ball of a specific color from the urn varies geometrically, with rate $q$. Then, on the stochastic model of a sequence of a specific number of random $q$-drawings of balls, a multivariate $q$-P\'{o}lya distribution is introduced and examined. Furthermore, in section 4, and on the stochastic model of a sequence of random $q$-drawings that is terminated with the occurrence of a specific number of balls of a given color, a multivariate inverse $q$-P\'{o}lya distribution is discussed.
\section{{\emph{q}}-Multinomial convolutions} \label{sec2}
\setcounter{equation}{0}
Let $x$ and $q$ be real numbers, with $q\neq 1$, and $r$ be an integer. The function of $x$, with parameter $q$, $[x]_q=(1-q^x)/(1-q)$ is called {\em $q$-number} and in particular $[r]_q$ is called {\em $q$-integer}. Also, the product $[x]_{r,q}=[x]_q[x-1]_q\cdots[x-r+1]_q$, $r=1,2,\ldots\,$, defines the {$q$-factorial of $x$ of order $r$. In particular, $[r]_q!=[1]_q[2]_q\cdots[r]_q$ is the {\em $q$-factorial} of $r$. The notion of $q$-factorial is extended to zero  order by $[x]_{0,q}=1$ and to negative order by $[x]_{-r,q}=1/[x+r]_{r,q}$, $r=1,2,\ldots\,$.

The {\em $q$-multinomial coefficient} is defined by
\begin{eqnarray}\label{eq2.1}
\genfrac[]{0pt}{}{x}{r_1,r_2,\ldots,r_k}_q=\frac{[x]_{r_1+r_2+\cdots+r_k,q}}{[r_1]_q![r_2]_q!\cdots[r_k]_q!}, \ \ r_j=0,1,\ldots,
\ \ j=1,2,\ldots,k.
\end{eqnarray}
Notice that a $q^{-1}$-number is readily expressed into a $q$-number by $[x]_{q^{-1}}=q^{-x+1}[x]_q$. Consequently,
\[
[x]_{r,q^{-1}}=q^{-xr+\binom{r+1}{2}}[x]_{r,q}, \ \ [r]_{q^{-1}}!=q^{-\binom{r}{2}}[r]_q!.
\]
Furthermore, setting $s_j=\sum_{i=1}^jr_i$ and $m_j=\sum_{i=j}^kr_i$, for $j=1,2,\ldots,k$, and using the expression
\[
-xs_k+\binom{s_k+1}{2}+\sum_{j=1}^k\binom{r_j}{2}=-\sum_{j=1}^k r_j(x-m_j)=-\sum_{j=1}^k r_j(x-s_j),
\]
it follows that
\begin{align}\label{eq2.2}
\genfrac{[}{]}{0pt}{}{x}{r_1,r_2,\ldots,r_k}_{q^{-1}}
&=q^{-\sum_{j=1}^k r_j(x-m_j)}\genfrac{[}{]}{0pt}{}{x}{r_1,r_2,\ldots,r_k}_q\nonumber\\
&=q^{-\sum_{j=1}^k r_j(x-s_j)}\genfrac{[}{]}{0pt}{}{x}{r_1,r_2,\ldots,r_k}_q.
\end{align}
Therefore, a formula involving $q$-numbers, $q$-factorials, and $q$-multinomial coefficients in the base $q$, with $1<q<\infty$, can be converted, with respect to the base, into a similar formula in the base $p=q^{-1}$, with $0<p<1$. Two versions of a recurrence relation for the $q$-multinomial coefficients, useful in the sequel, are quoted here for easy reference.

The $q$-multinomial coefficient satisfies the recurrence relation
\begin{align}\label{eq2.3}
\genfrac{[}{]}{0pt}{}{x}{r_1,r_2,\ldots,r_k}_q&=\genfrac{[}{]}{0pt}{}{x-1}{r_1,r_2,\ldots,r_k}_q
+q^{x-m_1}\genfrac{[}{]}{0pt}{}{x-1}{r_1-1,r_2,\ldots,r_k}_q\nonumber \\
&+q^{x-m_2}\genfrac{[}{]}{0pt}{}{x-1}{r_1,r_2-1,\ldots,r_k}_q
+\cdots+q^{x-m_k}\genfrac{[}{]}{0pt}{}{x-1}{r_1,r_2,\ldots,r_k-1}_q,
\end{align}
and alternatively, the recurrence relation
\begin{align}\label{eq2.4}
\genfrac{[}{]}{0pt}{}{x}{r_1,r_2,\ldots,r_k}_q&=q^{s_k}\genfrac{[}{]}{0pt}{}{x-1}{r_1,r_2,\ldots,r_k}_q
+\genfrac{[}{]}{0pt}{}{x-1}{r_1-1,r_2,\ldots,r_k}_q\nonumber \\
&+q^{s_1}\genfrac{[}{]}{0pt}{}{x-1}{r_1,r_2-1,\ldots,r_k}_q
+\cdots+q^{s_{k-1}}\genfrac{[}{]}{0pt}{}{x-1}{r_1,r_2,\ldots,r_k-1}_q,
\end{align}
for $r_j=0,1,\ldots\,$ and $j=1,2,\ldots,k$, with $m_j=\sum_{i=j}^kr_i$ and $s_j=\sum_{i=1}^jr_i$.

Recurrence relations (\ref{eq2.3}) and (\ref{eq2.4}), by replacing the base $q$ by $q^{-1}$, and using the first and the second expression in (\ref{eq2.2}), respectively, are expressed as
\begin{align} \label{eq2.5}
\genfrac{[}{]}{0pt}{}{x}{r_1,r_2,\ldots,r_k}_q=\,&q^{m_1}\genfrac{[}{]}{0pt}{}{x-1}{r_1,r_2,\ldots,r_k}_q
+q^{m_2}\genfrac{[}{]}{0pt}{}{x-1}{r_1-1,r_2,\ldots,r_k}_q \nonumber\\
&+q^{m_3}\genfrac{[}{]}{0pt}{}{x-1}{r_1,r_2-1,\ldots,r_k}_q
+\cdots+\genfrac{[}{]}{0pt}{}{x-1}{r_1,r_2,\ldots,r_k-1}_q.
\end{align}
and
\begin{align} \label{eq2.6}
\genfrac{[}{]}{0pt}{}{x}{r_1,r_2,\ldots,r_k}_q&=\genfrac{[}{]}{0pt}{}{x-1}{r_1,r_2,\ldots,r_k}_q
+q^{x-s_1}\genfrac{[}{]}{0pt}{}{x-1}{r_1-1,r_2,\ldots,r_k}_q \nonumber\\
&+q^{x-s_2}\genfrac{[}{]}{0pt}{}{x-1}{r_1,r_2-1,\ldots,r_k}_q
+\cdots+q^{x-s_k}\genfrac{[}{]}{0pt}{}{x-1}{r_1,r_2,\ldots,r_k-1}_q,
\end{align}
respectively. It should be noticed that reversing the order of labeling the arguments of the $q$-multinomials, these expressions may be transformed to (\ref{eq2.4}) and (\ref{eq2.3}), respectively.

Two versions of a multivariate $q$-Vandermonde formula are derived in the next theorem.
\begin{thm}\label{thm2.1}
Let $n$ be a positive integer and let $x_j$, $j=1,2,\ldots,k+1$, and $q$ be real numbers, with $q\neq 1$. Then,
\begin{eqnarray}\label{eq2.7}
[x_1+x_2+\cdots+x_{k+1}]_{n,q}=\sum\genfrac{[}{]}{0pt}{}{n}{r_1,r_2,\ldots,r_k}_q
q^{\sum_{j=1}^{k}(n-s_j)(x_j-r_j)}\prod_{j=1}^{k+1}[x_j]_{r_j,q},
\end{eqnarray}
and, alternatively,
\begin{eqnarray}\label{eq2.7a}
[x_1\!+\!x_2\!+\!\cdots\!+\!x_{k+1}\!+\!n\!-\!1]_{n,q}\!=\!\sum\genfrac{[}{]}{0pt}{}{n}{r_1,r_2,\ldots,r_k}_q
q^{\sum_{j=1}^{k}r_jz_j}\prod_{j=1}^{k+1}[x_j\!+\!r_j\!-\!1]_{r_j,q}.
\end{eqnarray}
Also,
\begin{eqnarray}\label{eq2.8}
[x_1+x_2+\cdots+x_{k+1}]_{n,q}=\sum\genfrac{[}{]}{0pt}{}{n}{r_1,r_2,\ldots,r_k}_q
q^{\sum_{j=1}^{k}r_j(z_j-(n-s_j))}\prod_{j=1}^{k+1}[x_j]_{r_j,q},
\end{eqnarray}
and, alternatively,
\begin{eqnarray}\label{eq2.8a}
[x_1\!+\!x_2\!+\!\cdots\!+\!x_{k+1}\!+\!n\!-\!1]_{n,q}\!=\!\sum\genfrac{[}{]}{0pt}{}{n}{r_1,r_2,\ldots,r_k}_q
\!q^{\sum_{j=1}^{k}x_j(n\!-\!s_j)}\!\prod_{j=1}^{k+1}[x_j\!+\!r_j\!-\!1]_{r_j,q},
\end{eqnarray}
where $s_j=\sum_{i=1}^jr_i$, $z_j=\sum_{i=j+1}^{k+1}x_i$, $j=1,2,\ldots,k$, and $r_{k+1}=n-s_k$, and the summation, in all four sums, is extended over all $r_j=0,1,\ldots,n$, $j=1,2,\ldots,k$, with $\sum_{i=1}^kr_i\leq n$.
\end{thm}
{\bf Proof}. Consider the sequence of multiple sums
\[
s_n(x_1,x_2,\ldots,x_{k+1};q)=\sum\genfrac{[}{]}{0pt}{}{n}{r_1,r_2,\ldots,r_k}_q
q^{\sum_{j=1}^{k}(n-s_j)(x_j-r_j)}\prod_{j=1}^{k+1}[x_j]_{r_j,q},
\]
for $n=1,2\ldots\,$, with initial value
\[
s_1(x_1,x_2,\ldots,x_{k+1};q)=[x_1]_q+\sum_{j=2}^{k+1}q^{\sum_{i=1}^{j-1}x_i}[x_j]_q=[x_1+x_2+\cdots+x_{k+1}]_q.
\]
Using recurrence relation (\ref{eq2.6}), with $x=n$, the sequence may be expressed as
\begin{align*}
s_n(x_1,x_2,\ldots,x_{k+1};q)
&=\sum\genfrac{[}{]}{0pt}{}{n-1}{r_1,r_2,\ldots,r_k}_q q^{\sum_{j=1}^{k}(n-s_j)(x_j-r_j)}\prod_{j=1}^{k+1}[x_j]_{r_j,q}\\
&+\sum\genfrac{[}{]}{0pt}{}{n-1}{r_1-1,r_2,\ldots,r_k}_q q^{n-s_1+\sum_{j=1}^{k}(n-s_j)(x_j-r_j)}\prod_{j=1}^{k+1}[x_j]_{r_j,q}\\
+\cdots&+\sum\genfrac{[}{]}{0pt}{}{n-1}{r_1,r_2,\ldots,r_k-1}_q q^{n-s_k+\sum_{j=1}^{k}(n-s_j)(x_j-r_j)}\prod_{j=1}^{k+1}[x_j]_{r_j,q}.
\end{align*}
Replacing $r_j-1$ by $r_j$ in the $(j+1)$th multiple sum, for $j=1,2,\ldots,k$, and then executing the multiplications, summations, and cancelations in the exponents of $q$, we get the expression
\begin{align*}
s_n&(x_1,x_2,\ldots,x_{k+1};q)\\
&=\sum\genfrac{[}{]}{0pt}{}{n-1}{r_1,r_2,\ldots,r_k}_q q^{\sum_{j=1}^{k}(n-1-s_j)(x_j-r_j)+\sum_{j=1}^k(x_j-r_j)}[x_{k+1}-r_{k+1}+1]_q\prod_{j=1}^{k+1}[x_j]_{r_j,q}\\
&+\sum\genfrac{[}{]}{0pt}{}{n-1}{r_1,r_2,\ldots,r_k}_q q^{\sum_{j=1}^{k}(n-1-s_j)(x_j-r_j)}[x_1-r_1]_q\prod_{j=1}^{k+1}[x_j]_{r_j,q}+\cdots\\
&+\sum\genfrac{[}{]}{0pt}{}{n-1}{r_1,r_2,\ldots,r_k}_q q^{\sum_{j=1}^{k}(n-1-s_j)(x_j-r_j)+\sum_{j=1}^{k-1}(x_j-r_j)}[x_k-r_k]_q\prod_{j=1}^{k+1}[x_j]_{r_j,q}.
\end{align*}
which, since
\begin{align*}
[x_1-r_1]_q&+q^{x_1-r_1}[x_2-r_2]_q+\cdots+q^{\sum_{j=1}^{k-1}(x_j-r_j)}[x_k-r_k]_q\\
&+q^{\sum_{j=1}^{k}(x_j-r_j)}[x_{k+1}-r_{k+1}+1]_q=[x_1,x_2,\ldots,x_{k+1}-n+1]_q
\end{align*}
implies for the sequence $s_n(x_1,x_2,\ldots,x_{k+1};q)$, $n=1,2\ldots\,$, the first-order recurrence relation
\[
s_n(x_1,x_2,\ldots,x_{k+1};q)=[x_1,x_2,\ldots,x_{k+1}-n+1]_q s_{n-1}(x_1,x_2,\ldots,x_{k+1};q),
\]
for $n=1,2\ldots\,$, with initial condition $s_1(x_1,x_2,\ldots,x_{k+1};q)=[x_1+x_2+\cdots+x_{k+1}]_q$. Applying it successively, it follows that $s_n(x_1,x_2,\ldots,x_{k+1};q)=[x_1+x_2+\cdots+x_{k+1}]_{n,q}$, and so (\ref{eq2.7}) is shown.

Formula (\ref{eq2.8}) may be derived by following the steps of the derivation of (\ref{eq2.7}) and using recurrence relation (\ref{eq2.4}), with $x=n$, and the expression
\begin{align*}
[x_{k+1}-r_{k+1}+1]_q&+q^{(x_{k+1}-r_{k+1})+1}[x_k-r_k]_q+\cdots+q^{\sum_{j=3}^{k+1}(x_j-r_j)+1}[x_2-r_2]_q\\
&+q^{\sum_{j=2}^{k+1}(x_j-r_j)+1}[x_1-r_1]_q=[x_1,x_2,\ldots,x_{k+1}-n+1]_q.
\end{align*}

The alternative formulae (\ref{eq2.7a}) and (\ref{eq2.8a}) are deduced from (\ref{eq2.7}) and (\ref{eq2.8}), respectively, by replacing $x_j$ by $-x_j$, for $j=1,2,\ldots,k+1$ and using the relations
\[
[-x]_{r,q}=(-1)^xq^{-xr-\binom{r}{2}}[x+r-1]_{r,q}, \ \ \binom{n}{2}=\sum_{j=1}^{k+1}\binom{r_j}{2}+\sum_{j=1}^kr_j(n-s_j).
\]

Two versions of a multivariate $q$-Cauchy formula, which by virtue of
\[
\genfrac{[}{]}{0pt}{}{n}{r_1,r_2,\ldots,r_k}_q=\frac{[n]_q!}{[r_1]_q![r_2]_q!\cdots[r_k]_q![r_{k+1}]_q!},\ \
\genfrac{[}{]}{0pt}{}{x_j}{r_j}_q=\frac{[x_j]_{r_j,q}}{[r_j]_q!},
\]
constitute reformulations of the corresponding two versions of a multivariate $q$-Vander\-monde formula, are stated in the following corollary of Theorem 2.1.
\begin{cor}\label{cor2.1}
Let $n$ be a positive integer and let $x_j$, $j=1,2,\ldots,k+1$, and $q$ be real numbers, with $q\neq 1$. Then,
\begin{eqnarray}\label{eq2.9}
\genfrac{[}{]}{0pt}{}{x_1+x_2+\cdots+x_{k+1}}{n}_q=\sum q^{\sum_{j=1}^{k}(n-s_j)(x_j-r_j)}\prod_{j=1}^{k+1}\genfrac{[}{]}{0pt}{}{x_j}{r_j}_q
\end{eqnarray}
and, alternatively,
\begin{eqnarray}\label{eq2.9a}
\genfrac{[}{]}{0pt}{}{x_1+x_2+\cdots+x_{k+1}+n-1}{n}_q=\sum q^{\sum_{j=1}^{k}r_jz_j}\prod_{j=1}^{k+1}\genfrac{[}{]}{0pt}{}{x_j+r_j-1}{r_j}_q.
\end{eqnarray}
Also,
\begin{eqnarray}\label{eq2.10}
\genfrac{[}{]}{0pt}{}{x_1+x_2+\cdots+x_{k+1}}{n}_q=\sum q^{\sum_{j=1}^{k}r_j(z_j-(n-s_j))}\prod_{j=1}^{k+1}\genfrac{[}{]}{0pt}{}{x_j}{r_j}_q
\end{eqnarray}
and, alternatively,
\begin{eqnarray}\label{eq2.10a}
\genfrac{[}{]}{0pt}{}{x_1+x_2+\cdots+x_{k+1}+n-1}{n}_q
=\sum q^{\sum_{j=1}^{k}x_j(n\!-\!s_j)}\prod_{j=1}^{k+1}\genfrac{[}{]}{0pt}{}{x_j+r_j-1}{r_j}_q,
\end{eqnarray}
where $s_j=\sum_{i=1}^jr_i$, $z_j=\sum_{i=j+1}^{k+1}x_i$, $j=1,2,\ldots,k$, and $r_{k+1}=n-s_k$ and the summation, in all four sums, is extended over all $r_j=0,1,\ldots,n$, $j=1,2,\ldots,k$, with $\sum_{i=1}^kr_i\leq n$.
\end{cor}
\begin{rem}
Additional expressions of the multivariate $q$-Cauchy formulae. {\em The alter\-native expressions (\ref{eq2.9a}) and (\ref{eq2.10a}), which are useful in probability theory, may be rewritten as
\begin{eqnarray}\label{eq2.9b}
\genfrac{[}{]}{0pt}{}{r+k}{n+k}_q=\sum q^{\sum_{j=1}^{k}(r_j-x_j)(n-y_j+k-j+1)}\prod_{j=1}^{k+1}\genfrac{[}{]}{0pt}{}{r_j}{x_j}_q
\end{eqnarray}
and
\begin{eqnarray}\label{eq2.10b}
\genfrac{[}{]}{0pt}{}{r+k}{n+k}_q
=\sum q^{\sum_{j=1}^{k}(x_j+1)(n-s_j-r+y_j)}\prod_{j=1}^{k+1}\genfrac{[}{]}{0pt}{}{r_j}{x_j}_q,
\end{eqnarray}
respectively, where $s_j=\sum_{i=1}^jr_i$, $y_j=\sum_{i=1}^{j}x_i$, $j=1,2,\ldots,k$, $r_{k+1}=r-s_k$, $x_{k+1}=r-y_k$, and the summation, in both sums, is extended over all $r_j=x_j,x_j+1,\ldots,r$, $j=1,2,\ldots,k$, with $\sum_{i=1}^kr_i\leq r$. Indeed, replacing, the bound variable $r_j$ by $r_j-(x_j-1)$ and the constant $x_j$ by $x_j+1$, for all $j=1,2,\ldots,k+1$, formulae (\ref{eq2.9a}) and (\ref{eq2.10a}), after some algebra, are transformed to (\ref{eq2.9b}) and (\ref{eq2.10b}), respectively.
}
\end{rem}

Two versions of a multivariate inverse $q$-Vandermonde formula are derived in the following theorem.
\begin{thm}\label{thm2.2}
Let $n$ be a positive integer and let $x_j$, $j=1,2,\ldots,k+1$, and $q$ be real numbers, with $q\neq 1$. Then,
\begin{eqnarray}\label{eq2.11}
\frac{1}{[x_{k+1}]_{n,q}}\!=\!\sum\genfrac{[}{]}{0pt}{}{n+s_k-1}{r_1,r_2,\ldots,r_k}_q
q^{\sum_{j=1}^{k}(n\!+\!s_k\!-\!s_j)(x_j-r_j)}\frac{\prod_{j=1}^k[x_j]_{r_j,q}}{[x_1+x_2+\cdots+x_{k+1}]_{n+s_k,q}},
\end{eqnarray}
provided $|q^{-x_{k+1}}|<1$, and
\begin{eqnarray}\label{eq2.12}
\frac{1}{[x_{k+1}]_{n,q}}\!=\!\sum\genfrac{[}{]}{0pt}{}{n+s_k-1}{r_1,r_2,\ldots,r_k}_q
q^{\sum_{j=1}^{k}r_j(z_j\!-\!s_k\!+\!s_j\!-\!n\!+\!1)}\frac{\prod_{j=1}^k[x_j]_{r_j,q}}{[x_1+x_2+\cdots+x_{k+1}]_{n+s_k,q}},
\end{eqnarray}
provided $|q^{x_{k+1}}|<1$, where $s_j=\sum_{i=1}^jr_i$ and $z_j=\sum_{i=j+1}^{k+1}x_i$, for $j=1,2,\ldots,k$, and the summation, in both sums, is extended over all $r_j=0,1,\ldots\,$, $j=1,2,\ldots,k$.
\end{thm}
{\bf Proof}. According to an inverse $q$-Vandermonde formula (Charalambides (2016), p. 14), it holds true
\[
\frac{1}{[x_{k+1}]_{n,q}}=\sum_{r_k=0}^\infty\genfrac{[}{]}{0pt}{}{n+r_k-1}{r_k}_q
q^{n(x_k-r_k)}\frac{[x_k]_{r_k,q}}{[x_k+x_{k+1}]_{n+r_k,q}}.
\]
Similarly,
\[
\frac{1}{[x_k+x_{k+1}]_{n+r_k,q}}=\sum_{r_{k-1}=0}^\infty\genfrac{[}{]}{0pt}{}{n+r_k+r_{k-1}-1}{r_{k-1}}_q
\frac{q^{(n+r_k)(x_{k-1}-r_{k-1})}[x_{k-1}]_{r_{k-1},q}}{[x_{k-1}+x_k+x_{k+1}]_{n+r_k+r_{k-1},q}}
\]
and finally,
\[
\frac{1}{[x_2+x_3+\cdots+x_{k+1}]_{n+s_k-s_1,q}}=\sum_{r_1=0}^\infty\genfrac{[}{]}{0pt}{}{n+s_k-1}{r_1}_q
\frac{q^{(n+s_k-s_1)(x_1-r_1)}[x_1]_{r_1,q}}{[x_1+x_2+\cdots+x_{k+1}]_{n+s_k,q}}.
\]
Applying these $k$ expansions, one after the other in the inner sum of each step, and using the relation
\[
\genfrac{[}{]}{0pt}{}{n+r_k-1}{r_k}_q\genfrac{[}{]}{0pt}{}{n+r_k+r_{k-1}-1}{r_{k-1}}_q\genfrac{[}{]}{0pt}{}{n+s_k-1}{r_1}_q
=\genfrac{[}{]}{0pt}{}{n+s_k-1}{r_1,r_2,\ldots,r_k}_q,
\]
expansion (\ref{eq2.11}) is obtained. The alternative expansion (\ref{eq2.12}), is similarly deduced by using the following inverse $q$-Vandermonde expansions (Charalambides (2016), p. 14)
\[
\frac{1}{[x_{j+1}+\cdots+x_{k+1}]_{n+s_k-s_j,q}}=\sum_{r_j=0}^\infty\genfrac{[}{]}{0pt}{}{n+s_k-s_{j-1}-1}{r_j}_q
\frac{q^{r_j(z_j-s_k+s_j-n+1)}[x_j]_{r_j,q}}{[x_j+\cdots+x_{k+1}]_{n+s_k-s_{j-1},q}},
\]
for $j=1,2,\ldots,k$, with $s_0=0$.
\section{Multivariate {\emph{q}}-P\'{o}lya distribution}\label{sec3}
\setcounter{equation}{0}
A multiple $q$-P\'{o}lya urn model may be introduced, by first defining a $q$-analogue of the notion of a random drawing of a ball from an urn.

Consider an urn containing $r$ balls, $\{b_1,b_2,\ldots,b_r\}$, of $k+1$ different ordered colors, with $r_\nu$ distinct balls of color $c_\nu$, $\{b_{s_{\nu-1}+1},b_{s_{\nu-1}+2},\ldots,b_{s_\nu}\}$, for $\nu=1,2,\ldots,k+1$, where $s_0=0$, $s_\nu=\sum_{i=1}^\nu r_i$, for $\nu=1,2,\ldots,k+1$, with $s_{k+1}=r$. A {\em random $q$-drawing (or $q$-selection)} of a ball from the urn is carried out as follows. Assume that the balls in the urn are forced to pass through a random mechanism, one by one, in the order $(b_1,b_2,\ldots,b_r)$ or in the reverse order $(b_r,b_{r-1},\ldots,b_1)$. Also, suppose that each passing ball may or may not be caught by the mechanism, with probabilities $p=1-q$ and $q$, respectively. The first caught ball is drawn out of the urn. In the case all balls in the urn pass through the mechanism and no ball is caught, the ball passing procedure is repeated, with the same order. Clearly, the probability that ball $b_x$ is drawn from the urn is given by
\[
\sum^\infty_{k=0}(1-q)q^{(x-1)+rk}=(1-q)q^{x-1}\sum^\infty_{k=0}q^{rk}=\frac{q^{x-1}}{[r]_q},
\]
or by
\[
\sum^\infty_{k=0}(1-q)q^{(r-x)+rk}=\frac{q^{r-x}}{[r]_q}=\frac{q^{-(x-1)}}{[r]_{q^{-1}}},
\]
where $0<q<1$, according to whether the ball passing order is $(b_1,b_2,\ldots,b_r)$ or $(b_r,b_{r-1},\ldots,b_1)$. Consequently, the probability function of the number $N_r$ on the drawn ball is given by
\[
p_{r}(x;q)=P(N_r=x)=\frac{q^{x-1}}{[r]_q},\ \ x=1,2,\ldots,r,
\]
where $0<q<1$ or $1<q<\infty$. Note that this is the probability function of the discrete $q$-uniform distribution on the set $\{1,2,\ldots,r\}$. Also, the probability $P_r(r_\nu;q)$, that a ball of color $c_\nu$ is drawn from the urn is given by
\[
P_r(r_\nu;q)=P(s_{\nu-1}<N_r\leq s_\nu)=\frac{q^{s_{\nu-1}}[r_\nu]_q}{[r]_q}=\frac{q^{-(r-s_\nu)}[r_\nu]_{q^{-1}}}{[r]_{q^{-1}}},
\]
for $\nu=1,2,\ldots,k+1$, with $s_0=0$, where $0<q<1$ or $1<q<\infty$. As expected, the sum of these probabilities, on using successively the relation $[s]_q+q^s[r]_q=[s+r]_q$, is obtained as
\[
\sum_{\nu=1}^{k+1}P_r(r_\nu;q)=\frac{1}{[r]_q}\sum_{\nu=1}^{k+1}q^{s_{\nu-1}}[r_\nu]_q=\frac{[r_1+r_2+\cdots+r_{k+1}]_q}{[r]_q}=1,
\]
where $0<q<1$ or $1<q<\infty$. Finally, notice that a random $q$-drawing of a ball, for $q\rightarrow 1$ and since
\[
\lim_{q\rightarrow 1}P_{r}(r_\nu;q)=\frac{r_\nu}{r},\ \  \nu=1,2,\ldots,k+1,
\]
reduces to the usual random drawing of a ball from the urn.

Furthermore, assume that random $q$-drawings of balls are sequentially carried out, one after the other, from an urn, initially containing $r$ balls of $k+1$ different colors, with $r_\nu$ distinct balls of color $c_\nu$, for $\nu=1,2,\ldots,k+1$, according to the following scheme. After each $q$-drawing, the drawn ball is placed back in the urn together with $m$ balls of the same color. Then, the conditional probability of drawing a ball of color $c_\nu$ at the $i$th $q$-drawing, given that $j_\nu-1$ balls of color $c_\nu$ and a total of $i_{\nu-1}$ balls of colors $c_1,c_2,\ldots,c_{\nu-1}$ are drawn in the previous $i-1$ $q$-drawings, is given by
\begin{align}\label{eq3.1}
p_{i,j_\nu}(i_{\nu-1})&=\frac{q^{s_{\nu-1}+mi_{\nu-1}}\big(1-q^{r_\nu+m(j_\nu-1)}\big)}{1-q^{r+m(i-1)}}
=\frac{q^{-m(\beta_{\nu-1}-i_{\nu-1})}[\alpha_\nu-j_\nu+1]_{q^{-m}}}{[\alpha-i+1]_{q^{-m}}},
\end{align}
for $j_\nu=1,2,\ldots,i$, $i_\nu=0,1,\ldots,i-1$, $i=1,2,\ldots\,$, and $\nu=1,2,\ldots,k+1$, with $i_0=0$, where $0<q<1$ or $1<q<\infty$ and $\alpha=-r/m$, $\alpha_\nu=-r_\nu/m$, $\beta_\nu=-s_\nu/m$, $\nu=1,2,\ldots,k+1$, $\beta_0=0$. Note that $i_\nu=i_{\nu-1}+j_\nu=j_1+j_2+\cdots+j_\nu$, for $\nu=1,2,\ldots,k+1$. This model, which for $q\rightarrow 1$ and since
\[
p_{i,j_\nu}=\lim_{q\rightarrow 1}p_{i,j_\nu}(i_{\nu-1})=\frac{r_\nu+m(j_\nu-1)}{r+m(i-1)}
=\frac{\alpha_\nu-j_\nu+1}{\alpha-i+1},\ \ \alpha_\nu=-\frac{r_\nu}{m}, \ \ \alpha=-\frac{r}{m},
\]
for $j_\nu=1,2,\ldots,i$, $i=1,2,\ldots\,$, and $\nu=1,2,\ldots,k+1$, reduces to the (classical) multiple P\'{o}lya urn model, may be called {\em multiple $q$-P\'{o}lya urn model}.
\begin{Def}\label{def3.1}
Let $X_\nu$ be the number of balls of color $c_\nu$ drawn in $n$ $q$-drawings in a multiple $q$-P\'{o}lya urn model, with conditional probability of drawing a ball of color $c_\nu$ at the $i$th $q$-drawing, given that $j_\nu-1$ balls of color $c_\nu$ and a total of $i_{\nu-1}$ balls of colors $c_1,c_2,\ldots,c_{\nu-1}$ are drawn in the previous $i-1$ $q$-drawings, is given by (\ref{eq3.1}), for $\nu=1,2,\ldots,k$. The distribution of the random vector $(X_1,X_2,\ldots, X_k)$ is called $k$-variate $q$-P\'{o}lya distribution, with parameters $n$, $(\alpha_1,\alpha_2,\ldots,\alpha_k)$, $\alpha$, and $q$.
\end{Def}

The probability function of the $k$-variate $q$-P\'{o}lya distribution is obtained in the following theorem.
\begin{thm}\label{thm3.1}
The probability function of the $k$-variate $q$-P\'{o}lya distribution, with parameters $n$, $(\alpha_1,\alpha_2,\ldots,\alpha_k)$, $\alpha$, and $q$, is given by
\begin{align}\label{eq3.2}
P(X_1\!=\!x_1,X_2\!=\,&x_2,\ldots,X_k\!=\!x_k)=q^{-m\sum_{j=1}^{k}(n-y_j)(\alpha_j-x_j)}\prod_{j=1}^{k+1}\genfrac{[}{]}{0pt}{}{\alpha_j}{x_j}_{q^{-m}}
\bigg/\genfrac{[}{]}{0pt}{}{\alpha}{n}_{q^{-m}}\nonumber\\
&=\genfrac{[}{]}{0pt}{}{n}{x_1,x_2,\ldots, x_k}_{q^{-m}}
q^{-m\sum_{j=1}^{k}(n-y_j)(\alpha_j-x_j)}\frac{\prod_{j=1}^{k+1}[\alpha_j]_{x_j,q^{-m}}}{[\alpha]_{n,q^{-m}}},
\end{align}
for $x_j=0,1,\ldots,n$, $j=1,2,\ldots,k$, with $\sum_{j=1}^k x_j\leq n$, and $0<q<1$ or $1<q<\infty$, where $x_{k+1}=n-\sum_{j=1}^k x_j$, $\alpha_{k+1}=\alpha-\sum_{j=1}^k \alpha_j$, and $y_j=\sum_{i=1}^{j}x_i$, for $j=1,2,\ldots,k$.
\end{thm}
{\bf Proof}. The probability function $p_n(x_1,x_2,\ldots,x_k)=P(X_1\!=\!x_1,X_2\!=\!x_2,\ldots,X_k\!=\!x_k)$, on using the total probability theorem, satisfies the recurrence relation
\begin{align*}
p_n(x_1,x_2,\ldots,x_k)=&\,p_{n-1}(x_1,x_2,\ldots,x_k)\frac{q^{-m(\beta_k-y_k)}[\alpha_{k+1}-x_{k+1}+1]_{q^{-m}}}{[\alpha-n+1]_{q^{-m}}}\\
&+p_{n-1}(x_1-1,x_2,\ldots,x_k)\frac{[\alpha_1-x_1+1]_{q^{-m}}}{[\alpha-n+1]_{q^{-m}}}\\
&+p_{n-1}(x_1,x_2-1,\ldots,x_k)\frac{q^{-m(\beta_1-y_1)}[\alpha_2-x_2+1]_{q^{-m}}}{[\alpha-n+1]_{q^{-m}}}+\cdots\\
&+p_{n-1}(x_1,x_2,\ldots,x_k-1)\frac{q^{-m(\beta_{k-1}-y_{k-1})}[\alpha_k-x_k+1]_{q^{-m}}}{[\alpha-n+1]_{q^{-m}}},
\end{align*}
for $x_j=1,2,\ldots,n$, $j=1,2,\ldots,k$ and $n=1,2,\ldots\,$, with $\sum_{j=1}^kx_j\leq n$, $x_{k+1}=n-\sum_{j=1}^kx_j$, and
$p_0(0,0,\ldots,0)=1$. Also,
\[
p_n(0,0,\ldots,0)=\frac{\prod_{i=1}^nq^{-m\beta_k}[a_{k+1}-i+1]_{q^{-m}}}{\prod_{i=1}^n[a-i+1]_{q^{-m}}}
=\frac{q^{-mn\beta_k}[a_{k+1}]_{n,q^{-m}}}{[a]_{n,q^{-m}}}.
\]
Clearly, the sequence
\begin{eqnarray}\label{eq3.3}
c_n(x_1,x_2,\ldots,x_k)
=q^{m\sum_{j=1}^{k}(n-y_j)(\alpha_j-x_j)}\frac{[\alpha]_{n,q^{-m}}}{\prod_{j=1}^{k+1}[\alpha_j]_{x_j,q^{-m}}}p_n(x_1,x_2,\ldots,x_k)
\end{eqnarray}
satisfies the recurrence relation
\begin{align*}
c_n(x_1,x_2,\ldots,x_k)=c_{n-1}(x_1,x_2,\ldots,x_k)&+q^{-m(n-y_1)}c_{n-1}(x_1-1,x_2,\ldots,x_k)\\
&+q^{-m(n-y_2)}c_{n-1}(x_1,x_2-1,\ldots,x_k)\\
+\cdots&+q^{-m(n-y_k)}c_{n-1}(x_1,x_2,\ldots,x_k-1),
\end{align*}
for $x_j=1,2,\ldots,n$, $j=1,2,\ldots,k$ and $n=1,2,\ldots\,$, with $\sum_{j=1}^kx_j\leq n$, and $c_0(0,0,\ldots,0)\linebreak=1$. Since this recurrence relation, according to (\ref{eq2.6}), uniquely determines the $q$-multinomial coefficient,
\[
c_n(x_1,x_2,\ldots,x_k)=\genfrac{[}{]}{0pt}{}{n}{x_1,x_2,\ldots, x_k}_{q^{-m}},
\]
the second part of expression (\ref{eq3.2}) is readily deduced from (\ref{eq3.3}). The first part of (\ref{eq3.2}), which is a reformulation of the second, is deduced by using the expressions
\[
\genfrac{[}{]}{0pt}{}{n}{x_1,x_2,\ldots,x_k}_q=\frac{[n]_q!}{[x_1]_q![x_2]_q!\cdots[x_k]_q![x_{k+1}]_q!},\ \
\genfrac{[}{]}{0pt}{}{\alpha_j}{x_j}_q=\frac{[\alpha_j]_{x_j,q}}{[x_j]_q!}.
\]
Note that the multivariate $q$-Vandermonde formula (\ref{eq2.7}) and, equivalently, the multivariate $q$-Cauchy formula (\ref{eq2.9}), guarantees that the probabilities (\ref{eq3.2}) sum to unity.

Certain (and not any) marginal and conditional distributions of a $k$-variate $q$-P\'{o}lya distribution are derived in the next theorem.
\begin{thm}\label{thm3.2}
Suppose that the random vector $(X_1,X_2,\ldots,X_k)$ obeys a $k$-variate $q$-P\'{o}lya distribution, with parameters $n$, $(\alpha_1,\alpha_2,\ldots,\alpha_k)$, $\alpha$, and $q$. Then,

(a) the marginal distribution of the random vector $(X_1,X_2,\ldots,X_\nu)$ is a $\nu$-variate $q$-P\'{o}lya, with parameters $n$, $(\alpha_1,\alpha_2,\ldots,\alpha_\nu)$, $\alpha$, and $q$, for $\nu=1,2,\ldots,k$, and

(b) the conditional distribution of the random vector $(X_{\nu},X_{\nu+1},\ldots,X_{\nu+\kappa-1})$, given that $(X_1,X_2,\ldots,X_{\nu-1})=(x_1,x_2,\ldots,x_{\nu-1})$ is a $\kappa$-variate $q$-P\'{o}lya, with parameters $n$, $(\alpha_\nu,\alpha_{\nu+1},\ldots,\alpha_{\nu+\kappa-1})$, $\alpha-\alpha_1-\alpha_2-\cdots-\alpha_{\nu-1}$, and $q$, for $\kappa=1,2,\ldots,k-\nu+1$ and $\nu=1,2,\ldots,k$.
\end{thm}
{\bf Proof}. (a) Summing the probabilities (\ref{eq3.2}) for $x_j=0,1,\ldots,n-y_{\nu}$, $j=\nu+1,\linebreak\nu+2,\ldots,k$, with
$x_{\nu+1}+x_{\nu+2}+\cdots+x_k\leq n-y_\nu$, and using (\ref{eq2.9}), we get
\begin{align*}
P(X_1&=x_1,X_2=x_2,\ldots,X_{\nu}=x_{\nu})=q^{-m\sum_{j=1}^{\nu}(n-y_j)(\alpha_j-x_j)}
\prod_{j=1}^{\nu}\genfrac{[}{]}{0pt}{}{\alpha_j}{x_j}_{q^{-m}}\\
&\times\sum q^{-m\sum_{j=\nu+1}^{k}(n-y_j)(\alpha_j-x_j)}\prod_{j=\nu+1}^{k}\genfrac{[}{]}{0pt}{}{\alpha_j}{x_j}_{q^{-m}}
\genfrac{[}{]}{0pt}{}{\alpha-\alpha_1-\cdots-\alpha_k}{n-x_1-\cdots-x_k}_{q^{-m}}\bigg/\genfrac{[}{]}{0pt}{}{\alpha}{n}_{q^{-m}}\\
&\hspace{1.3cm}=q^{-m\sum_{j=1}^{\nu}(n-y_j)(\alpha_j-x_j)}\prod_{j=1}^{\nu}\genfrac{[}{]}{0pt}{}{\alpha_j}{x_j}_{q^{-m}}
\genfrac{[}{]}{0pt}{}{\alpha-\alpha_1-\cdots-\alpha_\nu}{n-x_1-\cdots-x_\nu}_{q^{-m}}
\bigg/\genfrac{[}{]}{0pt}{}{\alpha}{n}_{q^{-m}},
\end{align*}
which is the probability function of a $\nu$-variate $q$-P\'{o}lya distribution, with parameters $n$, $(\alpha_1,\alpha_2,\ldots,\alpha_\nu)$, $\alpha$, and $q$.

(b) The probability function of the conditional distribution of the random vector $(X_{\nu},X_{\nu+1},\ldots,X_{\nu+\kappa-1})$, given that $(X_1,X_2,\ldots,X_{\nu-1})=(x_1,x_2,\ldots,x_{\nu-1})$, using the expression
\begin{align*}
P(X_\nu=x_\nu,\ldots,X_{\nu+\kappa-1}&=x_{\nu+\kappa-1}|X_1=x_1,X_2=x_2,\ldots,X_{\nu-1}=x_{\nu-1})\\
&=\frac{P(X_1=x_1,X_2=x_2,\ldots,X_{\nu+\kappa-1}=x_{\nu+\kappa-1})}{P(X_1=x_1,X_2=x_2,\ldots,X_{\nu-1}=x_{\nu-1})},
\end{align*}
is obtained as
\begin{align*}
P(X_\nu&=x_\nu,\ldots,X_{\nu+\kappa-1}=x_{\nu+\kappa-1}|X_1=x_1,X_2=x_2,\ldots,X_{\nu-1}=x_{\nu-1})\\
&=\frac{q^{-m\sum_{j=1}^{\nu+\kappa-1}(n-y_j)(\alpha_j-x_j)}\prod_{j=1}^{\nu+\kappa-1}\genfrac{[}{]}{0pt}{}{\alpha_j}{x_j}_{q^{-m}}
\genfrac{[}{]}{0pt}{}{\alpha-\alpha_1-\cdots-\alpha_{\nu+\kappa-1}}{n-x_1-\cdots-x_{\nu+\kappa-1}}_{q^{-m}}
\big/\genfrac{[}{]}{0pt}{}{\alpha}{n}_{q^{-m}}}{q^{-m\sum_{j=1}^{\nu-1}(n-y_j)(\alpha_j-x_j)}\prod_{j=1}^{\nu-1}
\genfrac{[}{]}{0pt}{}{\alpha_j}{x_j}_{q^{-m}}\genfrac{[}{]}{0pt}{}{\alpha-\alpha_1-\cdots-\alpha_{\nu-1}}{n-x_1-\cdots-x_{\nu-1}}_{q^{-m}}
\big/\genfrac{[}{]}{0pt}{}{\alpha}{n}_{q^{-m}}}\\
&=q^{-m\sum_{j=\nu}^{\nu+\kappa-1}(n-y_j)(\alpha_j-x_j)}\prod_{j=\nu}^{\nu+\kappa-1}\genfrac{[}{]}{0pt}{}{\alpha_j}{x_j}_{q^{-m}}
\genfrac{[}{]}{0pt}{}{\alpha-\alpha_\nu-\cdots-\alpha_{\nu+\kappa-1}}{n-x_\nu-\cdots-x_{\nu+\kappa-1}}_{q^{-m}}
\bigg/\genfrac{[}{]}{0pt}{}{\alpha}{n}_{q^{-m}},
\end{align*}
which is the probability function of a $\kappa$-variate $q$-P\'{o}lya distribution, with parameters $n$, $(\alpha_\nu,\alpha_{\nu+1},\ldots,\alpha_{\nu+\kappa-1})$, $\alpha-\alpha_1-\alpha_2-\cdots-\alpha_{\nu-1}$, and $q$.

The multivariate $q$-P\'{o}lya distribution, for large $r$, can be approximated by a $q$-multinomial distribution of the second kind, which is introduced and studied in Charalambides (2020). Specifically, the following limiting theorem is derived.
\begin{thm}\label{thm3.3}
Consider the multivariate $q$-P\'{o}lya distribution, with probability function $p_{n}(r;m,q)=P(X_1=x_1,X_2=x_2,\ldots,X_k=x_k)$ given by (\ref{eq3.2}).

For $0<q<1$, assume that
\begin{eqnarray}\label{eq3.4}
\lim_{r\rightarrow\infty}\frac{[r-s_j]_{q{-1}}}{[r-s_{j-1}]_{q{-1}}}=\theta_j, \ \ s_j=\sum_{i=1}^jr_i, \ \ j=1,2,\ldots,k,\ \ s_0=0,
\end{eqnarray}
and in the case of a negative integer $m$ assume, in addition, that $\theta_j<q^{-m(\nu-1)}$, $j=1,2,\ldots,k$, for some positive integer $\nu$. Then,
\begin{eqnarray}\label{eq3.5}
\lim_{r\rightarrow\infty}p_{n}(r;m,q)=\genfrac[]{0pt}{}{n}{x_1,x_2,\ldots, x_k}_{q^{m}}
\prod_{j=1}^k\theta_j^{n-y_j}\prod_{i_j=1}^{x_j}(1-\theta_jq^{m(i_j-1)}),
\end{eqnarray}
for $x_j=0,1,\ldots,n$, $j=1,2,\ldots,k$, with $x_1+x_2+\cdots+x_k\leq n$, where $y_j=\sum_{i=1}^{j}x_i$, $0<q<1$  and $0<\theta_j<1$, $j=1,2,\ldots,k$, in the case $m$ is a positive integer, or $0<\theta_j<q^{-m(\nu-1)}$, $j=1,2,\ldots,k$, for some positive integer $\nu\ge n$, in the case $m$ is a negative integer.

Also, for $1<q<\infty$, assume that
\begin{eqnarray}\label{eq3.6}
\lim_{r\rightarrow\infty}\frac{[r_j]_{q}}{[r-s_{j-1}]_{q}}=\lambda_j,  \ \ s_j=\sum_{i=1}^jr_i, \ \ j=1,2,\ldots,k,\ \ s_0=0,
\end{eqnarray}
and in the case of a negative integer $m$ assume, in addition, that $\lambda_j<q^{m(\nu-1)}$, $j=1,2,\ldots,k$, for some positive integer $\nu$. Then,
\begin{eqnarray}\label{eq3.7}
\lim_{r\rightarrow\infty}p_{n}(r;m,q)=\genfrac[]{0pt}{}{n}{x_1,x_2,\ldots, x_k}_{q^{-m}}
\prod_{j=1}^k\lambda_j^{x_j}\prod_{i_j=1}^{n-y_j}(1-\lambda_jq^{-m(i_j-1)}),
\end{eqnarray}
for $x_j=0,1,\ldots,n$, $j=1,2,\ldots,k$, with $x_1+x_2+\cdots+x_k\leq n$, where $y_j=\sum_{i=1}^{j}x_i$, $1<q<\infty$  and $0<\lambda_j<1$, in the case $m$ is a positive integer, or $0<\lambda_j<q^{m(\nu-1)}$, for some positive integer $\nu\ge n$, in the case $m$ is a negative integer.
\end{thm}
{\bf Proof}. For $0<q<1$, the probability function (\ref{eq3.2}),
\begin{align*}
p_{n}(r;m,q)=\genfrac{[}{]}{0pt}{}{n}{x_1,x_2,\ldots, x_k}_{q^{-m}}
\prod_{j=1}^{k}q^{-m(n-y_j)(\alpha_j-x_j)}
\frac{[\alpha_j]_{x_j,q^{-m}}[\alpha-\beta_j]_{n-y_j,q^{-m}}}{[\alpha-\beta_{j-1}]_{n-y_{j-1},q^{-m}}},
\end{align*}
using (\ref{eq2.2}), may be written as
\begin{align*}
p_{n}(r;m,q)&=\genfrac{[}{]}{0pt}{}{n}{x_1,x_2,\ldots, x_k}_{q^{m}}
\prod_{j=1}^{k}\prod_{i_j=1}^{x_j}(q^{-r_j-m(i_j-1)}-1)q^{-r+s_j+m(i_j-1)}\\
&\hspace{4.3cm}\frac{\prod_{i_j=1}^{n-y_j}(q^{-r+s_j-m(i_j-1)}-1)q^{m(i_j-1)}}{\prod_{i_{j-1}=1}^{n-y_{j-1}}(q^{-r+s_{j-1}-m(i_{j-1}-1)}-1)q^{m(i-1)}}.
\end{align*}
Moreover, by the assumption (\ref{eq3.4}), it follows that
\begin{align*}
\lim_{r\rightarrow\infty}\frac{(q^{-r_j-m(i_j-1)}-1)q^{-r+s_j+m(i_j-1)}}{q^{-r+s_{j-1}}-1}=1&-q^{m(i_j-1)}
\lim_{r\rightarrow\infty}\frac{q^{-r+s_j}-1}{q^{-r+s_{j-1}}-1}\\
&+\lim_{r\rightarrow\infty}\frac{q^{m(i_j-1)}-1}{q^{-r+s_{j-1}}-1}=1-\theta_jq^{m(i_j-1)},
\end{align*}
\begin{align*}
\lim_{r\rightarrow\infty}\frac{(q^{-r+s_j-m(i_j-1)}-1)q^{m(i_j-1)}}{q^{-r+s_{j-1}}-1}&=
\lim_{r\rightarrow\infty}\frac{q^{-r+s_j}-1}{q^{-r+s_{j-1}}-1}
-\lim_{r\rightarrow\infty}\frac{q^{m(i_j-1)}-1}{q^{-r+s_{j-1}}-1}=\theta_j,
\end{align*}
and
\begin{align*}
\lim_{r\rightarrow\infty}\frac{(q^{-r+s_{j-1}-m(i_j-1)}-1)q^{m(i_j-1)}}{q^{-r+s_{j-1}}-1}&=
1-\lim_{r\rightarrow\infty}\frac{q^{m(i_j-1)}-1}{q^{-r+s_{j-1}}-1}=1,
\end{align*}
for $j=1,2,\ldots,k$. Thus, dividing both the numerator and denominator of the $j$th factor in the last expression of the probability function (\ref{eq3.2}) by $(q^{-r+s_{j-1}}-1)^{n-y_{j-1}}$ and taking the limits as $r\rightarrow\infty$, the limiting expression (\ref{eq3.5}) is readily deduced.

For $1<q<\infty$, the probability function (\ref{eq3.2}),
\begin{align*}
p_{n}(r;m,q)=\genfrac{[}{]}{0pt}{}{n}{x_1,x_2,\ldots, x_k}_{q^{-m}}
\prod_{j=1}^{k}q^{-m(n-y_j)(\alpha_j-x_j)}
\frac{[\alpha_j]_{x_j,q^{-m}}[\alpha-\beta_j]_{n-y_j,q^{-m}}}{[\alpha-\beta_{j-1}]_{n-y_{j-1},q^{-m}}},
\end{align*}
may be written as
\begin{align*}
p_{n}(r;m,q)&=\genfrac{[}{]}{0pt}{}{n}{x_1,x_2,\ldots, x_k}_{q^{-m}}
\prod_{j=1}^{k}\prod_{i_j=1}^{x_j}(1-q^{r_j+m(i_j-1)})q^{-m(i_j-1)}\\
&\hspace{4.1cm}\times\frac{\prod_{i_j=1}^{n-y_j}(1-q^{r-s_j+m(i_j-1)})q^{r_j-m(i_j-1)}}
{\prod_{i_{j-1}=1}^{n-y_{j-1}}(1-q^{r-s_{j-1}+m(i_{j-1})})q^{-m(i_{j-1}-1)}}.
\end{align*}
Moreover, by the assumption (\ref{eq3.6}), it follows that
\begin{align*}
\lim_{r\rightarrow\infty}\frac{(q^{r_j+m(i_j-1)}-1)q^{-m(i_j-1)}}{q^{r-s_{j-1}}-1}=\lim_{r\rightarrow\infty}\frac{q^{r_j}-1}{q^{r-s_{j-1}}-1}
-\lim_{r\rightarrow\infty}\frac{q^{-m(i_j-1)}-1}{q^{r-s_{j-1}}-1}=\lambda_j,
\end{align*}
\begin{align*}
\lim_{r\rightarrow\infty}\frac{(q^{r-s_j+m(i_j-1)}-1)q^{r_j-m(i_j-1)}}{q^{r-s_{j-1}}-1}=1
&-q^{-m(i_j-1)}\lim_{r\rightarrow\infty}\frac{q^{r_j}-1}{q^{r-s_{j-1}}-1}\\
&-\lim_{r\rightarrow\infty}\frac{q^{-m(i_j-1)}-1}{q^{r-s_{j-1}}-1}=1-\lambda_jq^{-m(i_j-1)},
\end{align*}
and
\begin{align*}
\lim_{r\rightarrow\infty}\frac{(q^{r-s_{j-1}+m(i_j-1)}-1)q^{-m(i_j-1)}}{q^{r-s_{j-1}}-1}&=
1-\lim_{r\rightarrow\infty}\frac{q^{-m(i_j-1)}-1}{q^{r-s_{j-1}}-1}=1,
\end{align*}
for $j=1,2,\ldots,k$. Thus, dividing both the numerator and denominator of the $j$th factor in the last expression of the probability function (\ref{eq3.2}) by $(q^{r-s_{j-1}}-1)^{n-y_{j-1}}$ and taking the limits as $r\rightarrow\infty$, the limiting expression (\ref{eq3.7}) is readily deduced.

The multiple $q$-P\'{o}lya urn model in the particular case $m=0$ reduces to $q$-drawings with replacement and the distribution (\ref{eq3.2}) reduces to the classical multinomial distribution with probability of success of the $\nu$th kind $p_\nu=q^{s_{\nu-1}}[r_\nu]_q/[r]_q$, $\nu=1,2,\ldots,k+1$.

Also, for $m=-1$, the case corresponds to $q$-drawings without replacement and the probability function (\ref{eq3.2}) reduces to a
\begin{align}\label{eq3.8}
P(X_1=x_1,X_2=x_2,\ldots,&\,X_k=x_k)=q^{\sum_{j=1}^{k}(n-y_j)(r_j-x_j)}\prod_{j=1}^{k+1}\genfrac{[}{]}{0pt}{}{r_j}{x_j}_{q}
\bigg/\genfrac{[}{]}{0pt}{}{r}{n}_{q}\nonumber\\
&=\genfrac{[}{]}{0pt}{}{n}{x_1,x_2,\ldots, x_k}_{q}
q^{\sum_{j=1}^{k}(n-y_j)(r_j-x_j)}\frac{\prod_{j=1}^{k+1}[r_j]_{x_j,q}}{[r]_{n,q}},
\end{align}
for $x_j=0,1,\ldots,n$, $j=1,2,\ldots,k$, with $\sum_{j=1}^k x_j\leq n$, and $0<q<1$ or $1<q<\infty$, where $x_{k+1}=n-\sum_{j=1}^k x_j$, $r_{k+1}=r-\sum_{j=1}^k r_j$, and $y_j=\sum_{i=1}^{j}x_i$, $j=1,2,\ldots,k$. The distribution with probability function (\ref{eq3.8}) may be called {\em multivariate $q$-hypergeometric distribution}.

Furthermore, for $m=1$, the case to $q$-drawings with replacement and addition of another ball of the same color. The particular probability function may be deduced from (\ref{eq3.2}) by setting $\alpha_j=-r_j$, for $j=1,2\ldots,k$, $\alpha_{k+1}=-r+s_k$ and $x_{k+1}=n-y_k$, as
\[
P(X_1=x_1,X_2=x_2,\ldots,X_k=x_k)=\prod_{j=1}^{k}q^{(n-y_j)(r_j+x_j)}\genfrac{[}{]}{0pt}{}{-r_j}{x_j}_{q^{-1}}
\genfrac{[}{]}{0pt}{}{-r+s_k}{n-y_k}_{q^{-1}}\bigg/\genfrac{[}{]}{0pt}{}{-r}{n}_{q^{-1}}.
\]
which, using the expression
\[
\genfrac{[}{]}{0pt}{}{-r+s_k}{n-y_k}_{q^{-1}}\bigg/\genfrac{[}{]}{0pt}{}{-r}{n}_{q^{-1}}
=\prod_{j=1}^{k}\genfrac{[}{]}{0pt}{}{-r+s_j}{n-y_j}_{q^{-1}}\bigg/\genfrac{[}{]}{0pt}{}{-r+s_{j-1}}{n-y_{j-1}}_{q^{-1}},
\]
with $y_0=0$, becomes
\[
P(X_1=x_1,\ldots,X_k=x_k)=\prod_{j=1}^{k}q^{(n-y_j)(r_j+x_j)}\genfrac{[}{]}{0pt}{}{-r_j}{x_j}_{q^{-1}}
\genfrac{[}{]}{0pt}{}{-r+s_j}{n-y_j}_{q^{-1}}\bigg/\genfrac{[}{]}{0pt}{}{-r+s_{j-1}}{n-y_{j-1}}_{q^{-1}}.
\]
Then, since
\[
\genfrac{[}{]}{0pt}{}{-r_j}{x_j}_{q^{-1}}=(-1)^{x_j}q^{x_j+\binom{x_j}{2}}\genfrac{[}{]}{0pt}{}{r_j+x_j-1}{x_j}_q,
\]
\[
\genfrac{[}{]}{0pt}{}{-r+s_j}{n-y_j}_{q^{-1}}=(-1)^{n-y_j}q^{n-y_j+\binom{n-y_j}{2}}\genfrac{[}{]}{0pt}{}{r-s_j+n-y_j-1}{n-y_j}_{q^{-1}},
\]
and
\[
\binom{n-y_{j-1}}{2}=\binom{x_j+(n-y_j)}{2}=\binom{x_j}{2}+\binom{n-y_j}{2}+x_j(n-y_j),
\]
it takes the form
\[
P(X_1=x_1,X_2=x_2,\ldots,X_k=x_k)=\prod_{j=1}^{k}q^{r_j(n-y_j)}\frac{\genfrac{[}{]}{0pt}{}{r_j+x_j-1}{x_j}_q
\genfrac{[}{]}{0pt}{}{r-s_j+n-y_j-1}{n-y_j}_q}{\genfrac{[}{]}{0pt}{}{r-s_{j-1}+n-y_{j-1}-1}{n-y_{j-1}}_q},
\]
which after cancelations, reduces to
\begin{align}\label{eq3.9}
P(X_1\!=x_1,X_2\!=x_2,&\ldots,X_k\!=x_k)\!=\!q^{\sum_{j=1}^kr_j(n-y_j)}\prod_{j=1}^{k+1}\genfrac{[}{]}{0pt}{}{r_j+x_j-1}{x_j}_q
\bigg/\genfrac{[}{]}{0pt}{}{r+n-1}{n}_q\nonumber\\
&=\genfrac{[}{]}{0pt}{}{n}{x_1,x_2,\ldots, x_k}_{q}
q^{\sum_{j=1}^{k}r_j(n-y_j)}\frac{\prod_{j=1}^{k+1}[r_j+x_j-1]_{x_j,q}}{[r+n-1]_{n,q}},
\end{align}
for $x_j=0,1,\ldots,n$, $j=1,2,\ldots,k$, with $\sum_{j=1}^k x_j\leq n$, and $0<q<1$ or $1<q<\infty$, where $x_{k+1}=n-\sum_{j=1}^k x_j$, $r_{k+1}=r-\sum_{j=1}^k r_j$, and $y_j=\sum_{i=1}^{j}x_i$, $j=1,2,\ldots,k$. The distribution with probability function (\ref{eq3.9}) may be called {\em multivariate negative $q$-hypergeometric distribution}. Note that the probabilities (\ref{eq3.9}), according to (\ref{eq2.8a}), sum to unity.
\begin{exm}\label{exm3.1}
Distribution of the numbers of errors in the chapters of a manuscript. {\em Consider a manuscript of $k+1$ chapters (sections, parts), with chapter $c_\nu$ containing $r_\nu$ typographical errors, $\{e_{s_{\nu-1}+1},e_{s_{\nu-1}+2},\ldots,e_{s_\nu}\}$, for $\nu=1,2,\ldots,k+1$, where $s_0=0$, $s_\nu=\sum_{i=1}^\nu r_i$, for $\nu=1,2,\ldots,k+1$, with $s_{k+1}=r$. Assume that a proofreader reads the manuscript and when he/she finds an error corrects it and starts reading the manuscript from the beginning. Also, the proofreader starts reading the manuscript from the beginning when he/she reaches its end. Assume that the probability of finding any particular error is $p=1-q$. Clearly, the probability of finding an error in chapter $c_\nu$ at the first scan is given by $q^{s_{\nu-1}}[r_\nu]_q/[r]_q$, for $\nu=1,2,\ldots,k+1$, with $s_0=0$. Then, the conditional probability of finding (and correcting) an error in chapter $c_\nu$ at the $i$th scan, given that $j_\nu-1$ errors of chapter $c_\nu$ and a total of $i_{\nu-1}$ errors of chapters $c_1,c_2,\ldots,c_{\nu-1}$ are found in the previous $i-1$ scans, is given by
\begin{align*}
p_{i,j_\nu}(i_{\nu-1})=\frac{q^{s_{\nu-1}-i_{\nu-1}}[r_\nu-j_\nu+1]_q}{[r-i+1]_q},
\end{align*}
for $j_\nu=1,2,\ldots,i$, $i_\nu=0,1,\ldots,i-1$, $i=1,2,\ldots\,$, and $\nu=1,2,\ldots,k+1$, with $s_0=0$ and $i_0=0$, where $0<q<1$. Clearly, the joint distributions of the numbers $X_\nu$ of errors of chapter $c_\nu$ found (and corrected) in $n$ scans is the multivariate $q$-hypergeometric distribution, with probability function (\ref{eq3.8}).}
\end{exm}
\begin{exm}\label{exm3.2}
Random $q$-selection from a finite population. {\em Consider a finite population of $r$ people, classified into $k+1$ classes $c_j$, $j=1,2,\ldots,k$, with an unknown number of people in each class. Suppose that a sample of $n$ people is randomly $q$-selected from this population, without replacement. Let $x_j$ be the number of people of class $c_j$, for $j=1,2,\ldots,k$, in the sample. We are interested in the probability that the number of people of the population who belong in class $c_j$ equals $r_j$, for $j=1,2,\ldots,k$.

Let $X_j$ and $R_j$ be the numbers of people of class $c_j$, for $j=1,2,\ldots,k$, in the sample and the population, respectively. The conditional distribution of the random vector ${\textbf{X}}=(X_1,X_2,\ldots,X_k)$, given that the random vector ${\textbf{R}}=(R_1,R_2,\ldots,R_k)$ equals ${\textbf{r}}=(r_1,r_2,\ldots,r_k)$, is the multivariate $q$-hypergeometric distribution, with probability function (\ref{eq3.8}),
\[
f_{{\textbf{X}}|{\textbf{R}}}(x_1,x_2,\ldots,x_k|r_1,r_2,\ldots,r_k)=q^{\sum_{j=1}^{k}(n-y_j)(r_j-x_j)}
\prod_{j=1}^{k+1}\genfrac{[}{]}{0pt}{}{r_j}{x_j}_{q}\bigg/\genfrac{[}{]}{0pt}{}{r}{n}_{q},
\]
for $x_j=0,1,\ldots,n$, $j=1,2,\ldots,k$, with $\sum_{j=1}^k x_j\leq n$, and $0<q<1$ or $1<q<\infty$, where $x_{k+1}=n-\sum_{j=1}^k x_j$, $r_{k+1}=r-\sum_{j=1}^k r_j$, and $y_j=\sum_{i=1}^{j}x_i$, $j=1,2,\ldots,k$. The required probability is given by the value of the conditional probability function of the random vector ${\textbf{R}}=(R_1,R_2,\ldots,R_k)$, given ${\textbf{X}}=(X_1,X_2,\ldots,X_k)$, at the point $\textbf{r}=(r_1,r_2,\ldots,r_k)$. This conditional probability function is given
\begin{align*}
f_{{\textbf{R}}|{\textbf{X}}}(r_1,r_2,\ldots,r_k|x_1,x_2,\ldots,x_k)
&=\frac{f_{{\textbf{R}},{\textbf{X}}}(r_1,r_2,\ldots,r_k,x_1,x_2,\ldots,x_k)}{f_{{\textbf{X}}}(x_1,x_2,\ldots,x_k)}\\
&=\frac{f_{{\textbf{R}}}(r_1,r_2,\ldots,r_k)f_{{\textbf{X}}|{\textbf{R}}}(x_1,x_2,\ldots,x_k|r_1,r_2,\ldots,r_k)}
{f_{{\textbf{X}}}(x_1,x_2,\ldots,x_k)}
\end{align*}
and the probability function of the random vector ${\textbf{X}}=(X_1,X_2,\ldots,X_k)$ is
\[
f_{{\textbf{X}}}(x_1,x_2,\ldots,x_k)
=\sum f_{{\textbf{R}}}(r_1,r_2,\ldots,r_k)f_{{\textbf{X}}|{\textbf{R}}}(x_1,x_2,\ldots,x_k|r_1,r_2,\ldots,r_k),
\]
where the summation is extended over all $r_j=x_j,x_j+1,\ldots,r$, $j=1,2,\ldots,k$, with $\sum_{j=1}^kr_j\le r$. Thus, for the calculation of the probability in question, the additional knowledge of the distribution of the random vector ${\textbf{R}}=(R_1,R_2,\ldots,R_k)$ is required. Assume that this distribution is the $k$-variate discrete $q$-uniform with probability function ({\em Bose-Einstein $q$-stochastic model} ({\em $q$-statistic}))
\[
f_{{\textbf{R}}}(r_1,r_2,\ldots,r_k)=\prod_{j=1}^{k+1}q^{(k-j+1)r_j}\bigg/\genfrac{[}{]}{0pt}{}{r+k}{k}_{q},
\]
for $r_j=0,1,\ldots,r$, $j=1,2,\ldots,k$, with $\sum_{j=1}^kr_j\leq r$, where $r_{k+1}=r-\sum_{j=1}^kr_j$. Hence
\[
f_{{\textbf{X}}}(x_1,x_2,\ldots,x_k)
=\sum q^{\sum_{j=1}^{k}(n-y_j)(r_j-x_j)+(k-j+1)r_j}
\frac{\prod_{j=1}^{k+1}\genfrac{[}{]}{0pt}{}{r_j}{x_j}_{q}}{\genfrac{[}{]}{0pt}{}{r+k}{k}_{q}\genfrac{[}{]}{0pt}{}{r}{n}_{q}},
\]
where the summation is extended over all $r_j=x_j,x_j+1,\ldots,r$, $j=1,2,\ldots,k$, with $\sum_{j=1}^kr_j\le r$, and $0<q<1$ or $1<q<\infty$, where $r_{k+1}=r-\sum_{j=1}^kr_j$, $x_{k+1}=n-\sum_{j=1}^kx_j$, and $y_j=\sum_{i=1}^jx_i$. Then, using the $q$-Cauchy formula,
\[
\sum q^{\sum_{j=1}^{k}(r_j-x_j)(n-y_j+k-j+1)}
\prod_{j=1}^{k+1}\genfrac{[}{]}{0pt}{}{r_j}{x_j}_{q}=\genfrac{[}{]}{0pt}{}{r+k}{n+k}_{q}
\]
where the summation is extended over all $r_j=x_j,x_j+1,\ldots,r$, $j=1,2,\ldots,k$, with $\sum_{j=1}^kr_j\le r$, and $0<q<1$ or $1<q<\infty$,  where $r_{k+1}=r-\sum_{j=1}^kr_j$, $x_{k+1}=n-\sum_{j=1}^kx_j$, and $y_j=\sum_{i=1}^jx_i$, the probability function of the random vector ${\textbf{X}}=(X_1,X_2,\ldots,X_k)$, is deduced as
\[
f_{{\textbf{X}}}(x_1,x_2,\ldots,x_k)=\prod_{j=1}^{k+1}q^{(k-j+1)x_j}\genfrac{[}{]}{0pt}{}{r+k}{n+k}_{q}
\bigg/\bigg(\genfrac{[}{]}{0pt}{}{r+k}{k}_{q}\genfrac{[}{]}{0pt}{}{r}{n}_{q}\bigg).
\]
Since
\[
\genfrac{[}{]}{0pt}{}{r+k}{k}_{q}\genfrac{[}{]}{0pt}{}{r}{n}_{q}=\genfrac{[}{]}{0pt}{}{r+k}{n+k}_{q}\genfrac{[}{]}{0pt}{}{n+k}{k}_{q},
\]
the last expression reduces to
\[
f_{{\textbf{X}}}(x_1,x_2,\ldots,x_k)=\prod_{j=1}^{k+1}q^{(k-j+1)x_j}\bigg/\genfrac{[}{]}{0pt}{}{n+k}{k}_{q},
\]
for $x_j=0,1,\ldots,n$, $j=1,2,\ldots,k$, with $\sum_{j=1}^kr_j\le r$, and $0<q<1$ or $1<q<\infty$, which is a $k$-variate discrete $q$-uniform probability function. Therefore, the required conditional probability function of the random vector ${\textbf{R}}=(R_1,R_2,\ldots,R_k)$, given that ${\textbf{X}}=(X_1,X_2,\ldots,X_k)$ is given by
\[
f_{{\textbf{R}}|{\textbf{X}}}(r_1,r_2,\ldots,r_k|x_1,x_2,\ldots,x_k)
=q^{\sum_{j=1}^{k}(n-y_j)(r_j-x_j)}
\prod_{j=1}^{k+1}\genfrac{[}{]}{0pt}{}{r_j}{x_j}_{q}\bigg/\genfrac{[}{]}{0pt}{}{r+k}{n+k}_{q},
\]
for $r_j=x_j,x_j+1,\ldots,r$, $j=1,2,\ldots,k$, with $\sum_{j=1}^kr_j\le r$, and $0<q<1$ or $1<q<\infty$,  where $r_{k+1}=r-\sum_{j=1}^kr_j$, $x_{k+1}=n-\sum_{j=1}^kx_j$, and $y_j=\sum_{i=1}^jx_i$.}
\end{exm}

The multivariate $q$-hypergeometric distribution may be obtained as the conditional distribution of $k$ independent $q$-binomial distributions of the first kind, given their sum with another $q$-binomial distribution of the first kind independent of them. Precisely, the following theorem is shown.

\begin{thm}\label{thm3.4}
Consider a sequence of independent Bernoulli trials and assume that the probability of success at the $i$th trial is given by
\[
p_i=\frac{\theta q^{i-1}}{1+\theta q^{i-1}}, \ \ i=1,2,\ldots,\ \ 0<q<1 \ \ \text{or}\ \ 1<q<\infty.
\]
Let $X_j$ be the number of successes after the $(s_{j-1})$th trial and until the $(s_j)$th trial, for $j=1,2,\ldots,k+1$, with $s_0=0$, $s_j=\sum_{i=1}^jr_i$, $j=1,2,\ldots,k+1$, and $s_{k+1}=r$. Then, the conditional probability function of the random vector $(X_1,X_2,\ldots,X_k)$, given that $X_1+X_2+\cdots+X_{k+1}=n$, is the multivariate $q$-hypergeometric distribution with probability function (\ref{eq3.8}).
\end{thm}
{\bf Proof}. The random  variables $X_j$, $j=1,2,\ldots,k+1$, are independent, with probability function, according to Theorem 2.1 in Charalambides (2016), is given by
\[
P(X_j=x_j)=\genfrac{[}{]}{0pt}{}{r_j}{x_j}_q\frac{(\theta q^{s_{j-1}})^{x_j}q^{\binom{x_j}{2}}}{\prod_{i=1}^{r_j}(1+\theta q^{s_{j-1}+i-1})},
\ \ x_j=0,1,\ldots,r_j, \ \ j=1,2,\ldots,k+1.
\]
Similarly, the probability function of the sum $Y=X_1+X_2+\cdots+X_{k+1}$, which is the number of successes in $r$ trials, is
\[
P(Y=n)=\genfrac{[}{]}{0pt}{}{r}{n}_q\frac{\theta^n q^{\binom{n}{2}}}{\prod_{i=1}^r(1+\theta q^{i-1})},\ \ n=0,1,\ldots,r.
\]
Then, the joint conditional probability function of the random vector $(X_1,X_2,\ldots,X_k)$, given that $Y=n$,
\[
P(X_1=x_1,\ldots,X_k=x_k|Y=n)=\frac{P(X_1=x_1)\cdots P(X_k=x_k)P(X_{k+1}=n-y_k)}{P(Y=n)},
\]
on using these expressions, is obtained as
\[
P(X_1=x_1,X_2=x_2,\ldots,X_k=x_k|Y=n)=q^{c_k}\prod_{j=1}^{k+1}\genfrac{[}{]}{0pt}{}{r_j}{x_j}_{q}
\bigg/\genfrac{[}{]}{0pt}{}{r}{n}_{q},
\]
where
\[
c_k=\sum_{i=1}^kx_is_{i-1}+(n-y_k)s_k+\sum_{j=1}^k\binom{x_j}{2}+\binom{n-y_k}{2}-\binom{n}{2}.
\]
Thus, after some algebraic manipulations, it reduces to
\begin{align*}
c_k&=n\sum_{j=1}^kr_j-\sum_{i=1}^kx_i(s_k-s_{i-1})+\sum_{j=1}^k\binom{x_j}{2}+\binom{y_k+1}{2}-ny_k\\
&=\sum_{j=1}^kr_j(n-y_j)-\sum_{j=1}^kx_j(n-y_j)=\sum_{j=1}^k(n-y_j)(y_j-x_j),
\end{align*}
and the derivation of (\ref{eq3.8}) is completed.

Furthermore, the multivariate negative $q$-hypergeometric distribution may be obtained as the conditional distribution of $k$ independent negative $q$-binomial distributions of the second kind, given their sum with another negative $q$-binomial distribution of the second kind independent of them, according to the following theorem.

\begin{thm}\label{thm3.5}
Consider a sequence of independent Bernoulli trials and assume that the conditional probability of success at a trial, given that $j-1$ successes occur in the previous trials, is given by
\[
p_j=1-\theta q^{j-1}, \ \ j=1,2,\ldots,\ \ 0<\theta<1, \ \ 0<q<1 \ \ \text{or}\ \ 1<q<\infty,
\]
where, for $1<q<\infty$, the number $j$ of successes is restricted by $j\leq m=-\log\theta/\log q$.
Let $W_j$ be the number of failures after the $(s_{j-1})$th success and until the occurrence of the $(s_j)$th success, for $j=1,2,\ldots,k+1$, with $s_0=0$, $s_j=\sum_{i=1}^jr_i$, $j=1,2,\ldots,k+1$, and $s_{k+1}=r$, where $r\le m$ in the case $1<q<\infty$. Then, the conditional probability function of the random vector $(W_1,W_2,\ldots,W_k)$, given that $W_1+W_2+\cdots+W_{k+1}=n$, is the multivariate negative $q$-hypergeometric distribution with probability function (\ref{eq3.9}).
\end{thm}
{\bf Proof}. The random  variables $W_j$, $j=1,2,\ldots,k+1$, are independent, with probability function, according to Theorem 3.1 in Charalambides (2016), is given by
\[
P(W_j=w_j)=\genfrac{[}{]}{0pt}{}{r_j+w_j-1}{w_j}_q(\theta q^{s_{j-1}})^{w_j}\prod_{i=1}^{r_j}(1-\theta q^{s_{j-1}+i-1}),
\ \ w_j=0,1,\ldots,
\]
for all $j=1,2,\ldots,k+1$. Similarly, the probability function of the sum $U=W_1+W_2+\cdots+W_{k+1}$, which is the number of failures until the occurrence of the $r$th success, is
\[
P(U=n)=\genfrac{[}{]}{0pt}{}{r+n-1}{n}_q\theta^n \prod_{i=1}^r(1-\theta q^{i-1}),\ \ n=0,1,\ldots\,.
\]
Then, the joint conditional probability function of the random vector $(W_1,W_2,\ldots,W_k)$, given that $U=n$,
\[
P(W_1=w_1,\ldots,W_k=w_k|U=n)=\frac{P(W_1=w_1)\cdots P(W_k=w_k)P(W_{k+1}=n-u_k)}{P(U=n)},
\]
on using these expressions, is obtained as
\[
P(W_1=w_1,\ldots,W_k=w_k|U=n)=q^{c_k}\prod_{j=1}^{k+1}\genfrac{[}{]}{0pt}{}{r_j+w_j-1}{w_j}_q
\bigg/\genfrac{[}{]}{0pt}{}{r+n-1}{n}_q,
\]
where
\[
c_k=\sum_{i=1}^kw_is_{i-1}+(n-u_k)s_k, \ \ u_j=\sum_{i=1}^jw_i,\ \ j=1,2,\ldots,k.
\]
Thus, after some algebra, it reduces to
\[
c_k=n\sum_{j=1}^kr_j-\sum_{i=1}^kw_i(s_k-s_{i-1})=n\sum_{j=1}^kr_j-\sum_{j=1}^kr_ju_j=\sum_{j=1}^kr_j(n-u_j)
\]
and the derivation of (\ref{eq3.9}) is completed.
\section{Multivariate inverse {\emph{q}}-P\'{o}lya distribution}\label{sec4}
\setcounter{equation}{0}
Consider again the multiple $q$-P\'{o}lya urn model. Specifically, assume that random $q$-drawings of balls are sequentially carried out,  one after the other, from an urn, initially containing $r_\nu$ balls of color $c_\nu$, for $\nu=1,2,\ldots,k+1$, according to the following scheme. After each $q$-drawing the drawn ball is placed back in the urn together with $k$ balls of the same color. Assume that the conditional probability of drawing a ball of color $c_\nu$ at the $i$th $q$-drawing, given that $j_\nu-1$ balls of color $c_\nu$ and a total of $i_{\nu-1}$ balls of colors $c_1,c_2,\ldots,c_{\nu-1}$ are drawn in the previous $i-1$ $q$-drawings, is given by (\ref{eq3.1}). In this section the interest is turned to the study of the particular numbers of balls of colors $c_1,c_2,\ldots,c_k$ drawn until the $n$th ball of color $c_{k+1}$ is drawn. For this reason, the following definition is introduced.
\begin{Def}\label{def4.1}
Let $W_\nu$ be the number of balls of color $c_\nu$ drawn until the $n$th ball of color $c_{k+1}$ is drawn in a multiple $q$-P\'{o}lya urn model, with the conditional probability of drawing a ball of color $c_\nu$ at the $i$th $q$-drawing, given that $j_\nu-1$ balls of color $c_\nu$ and a total of $i_{\nu-1}$ balls of colors $c_1,c_2,\ldots,c_{\nu-1}$ are drawn in the previous $i-1$ $q$-drawings, given by (\ref{eq3.1}), for $\nu=1,2,\ldots,k$. The distribution of the random vector
$(W_1, W_2,\ldots,W_k)$ is called multivariate inverse $q$-P\'{o}lya distribution, with parameters $n$, $(\alpha_1,\alpha_2,\ldots,\alpha_k)$, $\alpha$, and $q$.
\end{Def}

The probability function of the $k$-variate inverse $q$-P\'{o}lya distribution is obtained in the following theorem.
\begin{thm}\label{thm4.1}
The probability function of the $k$-variate inverse $q$-P\'{o}lya distribution, with parameters $n$, $(\alpha_1,\alpha_2,\ldots,\alpha_k)$, $\alpha$, and $q$, is given by
\begin{align}\label{eq4.1}
P&(W_1\!=\!w_1,W_2\!=w_2,\ldots,W_k\!=\!w_k)\nonumber\\
&=\genfrac{[}{]}{0pt}{}{n+u_k-1}{w_1,w_2,\ldots, w_k}_{q^{-m}}
q^{-m\sum_{j=1}^{k}(n+u_k-u_j)(\alpha_j-w_j)}\frac{\prod_{j=1}^k[\alpha_j]_{w_j,q^{-m}}[\alpha_{k+1}]_{n,q^{-m}}}
{[\alpha]_{n+u_k,q^{-m}}},
\end{align}
for $w_j=0,1,\ldots\,$, $j=1,2,\ldots,k$, and $0<q<1$ or $1<q<\infty$, where $\alpha_{k+1}=\alpha-\sum_{j=1}^k\alpha_j$, and $u_j=\sum_{i=1}^jw_i$, for $j=1,2,\ldots,k$.
\end{thm}
{\bf Proof}. The probability function of the $k$-variate inverse $q$-P\'{o}lya distribution is closely connected to the probability function $k$-variate $q$-P\'{o}lya distribution. Specifically,
\[
P(W_1=w_1,W_2=w_2,\ldots,W_k=w_k)=p_{n+u_k-1}(w_1,w_2,\ldots, w_k)p_{n+u_k,n},
\]
where $p_{n+u_k-1}(w_1,w_2,\ldots, w_k)$ is the probability of drawing $w_\nu$ balls of color $c_\nu$, for all $\nu=1,2,\ldots,k$, and $n-1$ balls of color $c_{k+1}$ in $n+u_k-1$ $q$-drawings and $p_{n+u_k,n}=q^{-m(\beta_k-u_k)}[a_{k+1}-n+1]_{q^{-m}}/[a-n-u_k+1]_{q^{-m}}$ is the conditional probability of drawing a ball of color $c_{k+1}$ at the $(n+u_k)$th $q$-drawing, given that $n-1$ balls of color $c_{k+1}$ and a total of $u_k$ balls of colors $c_1,c_2,\ldots,c_k$ are drawn in the previous $n+u_k-1$ $q$-drawings. Thus using (\ref{eq3.2}), expression (\ref{eq4.1}) is deduced. Note that the multivariate inverse $q$-Vandermonde formula (\ref{eq2.11}) guarantees that the probabilities (\ref{eq4.1}) sum to unity.

The multivariate inverse $q$-P\'{o}lya distribution, for large $r$, can be approximated by a negative $q$-multinomial distribution of the second kind, which is introduced and studied in Charalambides (2020). Specifically, the following limiting theorem is derived.
\begin{thm}\label{thm5.2}
Consider the multivariate inverse $q$-P\'{o}lya distribution, with probability function $q_{n}(r;m,q)=P(W_1=w_1,W_2=w_2,\ldots,W_k=w_k)$ given by (\ref{eq4.1}).

For $0<q<1$, assume that the limiting expression (\ref{eq3.4}) holds true. Then,
\begin{eqnarray}\label{eq4.2}
\lim_{r\rightarrow\infty}q_{n}(r;m,q)=\genfrac[]{0pt}{}{n+u_k-1}{w_1,w_2,\ldots, w_k}_{q^{m}}
\prod_{j=1}^k\theta_j^{n+u_k-u_j}q^{mw_j}\prod_{i_j=1}^{w_j}(1-\theta_jq^{m(i_j-1)}),
\end{eqnarray}
for $w_j=0,1,\ldots\,$, and $j=1,2,\ldots,k$, where $u_j=\sum_{i=1}^jw_i$, $0<q<1$ and $0<\theta_j<1$, $j=1,2,\ldots,k$, in the case $m$ is a positive integer, or $0<\theta_j<q^{-m(\nu-1)}$, $j=1,2,\ldots,k$, for some positive integer $\nu\ge n$, in the case $m$ is a negative integer.

Also, for $1<q<\infty$, assume that the limiting expression (\ref{eq3.6}) holds true. Then,
\begin{eqnarray}\label{eq4.3}
\lim_{r\rightarrow\infty}q_{n}(r;m,q)=\genfrac[]{0pt}{}{n+u_k-1}{w_1,w_2,\ldots, w_k}_{q^{-m}}
\prod_{j=1}^k\lambda_j^{w_j}\prod_{i_j=1}^{n+u_k-u_j}(1-\lambda_jq^{-m(i_j-1)}),
\end{eqnarray}
for $w_j=0,1,\ldots\,$, and $j=1,2,\ldots,k$, where $u_j=\sum_{i=1}^{j}w_i$, $1<q<\infty$  and $0<\lambda_j<1$, $j=1,2,\ldots,k$, in the case $m$ is a positive integer, or $0<\lambda_j<q^{m(\nu-1)}$, $j=1,2,\ldots,k$, for some positive integer $\nu\ge n$, in the case $m$ is a negative integer.
\end{thm}
{\bf Proof}. The probability function (\ref{eq4.1}), using (\ref{eq2.2}), may be written as
\begin{align*}
q_{n}(r;m,q)&=\genfrac{[}{]}{0pt}{}{n+u_k-1}{w_1,w_2,\ldots, w_k}_{q^{-m}}
\prod_{j=1}^{k}q^{-m(n+u_k-u_j)(\alpha_j-w_j)}[\alpha_j]_{w_j,q^{-m}}\\
&\hspace{4.6cm}\times\frac{[\alpha-\beta_j]_{n+u_k-u_j,q^{-m}}}{[\alpha-\beta_{j-1}]_{n+u_k-u_{j-1},q^{-m}}}\\
&=\genfrac{[}{]}{0pt}{}{n+u_k-1}{w_1,w_2,\ldots, w_k}_{q^{m}}
\prod_{j=1}^{k}q^{mw_j}\prod_{i_j=1}^{w_j}(q^{-r_j-m(i_j-1)}-1)q^{-r+s_j+m(i_j-1)}\\
&\hspace{4.3cm}\times\frac{\prod_{i_j=1}^{n+u_k-u_j}(q^{-r+s_j-m(i_j-1)}-1)q^{m(i_j-1)}}
{\prod_{i_{j-1}=1}^{n+u_k-u_{j-1}}(q^{-r+s_{j-1}-m(i_{j-1}-1)}-1)q^{m(i-1)}}.
\end{align*}
and, alternatively, as
\begin{align*}
q_{n}(r;m,q)&=\genfrac{[}{]}{0pt}{}{n+u_k-1}{w_1,w_2,\ldots, w_k}_{q^{-m}}
\prod_{j=1}^{k}q^{-m(n+u_k-u_j)(\alpha_j-w_j)}[\alpha_j]_{w_j,q^{-m}}\\
&\hspace{4.6cm}\times\frac{[\alpha-\beta_j]_{n+u_k-u_j,q^{-m}}}{[\alpha-\beta_{j-1}]_{n+u_k-u_{j-1},q^{-m}}}\\
&=\genfrac{[}{]}{0pt}{}{n+u_k-1}{w_1,w_2,\ldots, w_k}_{q^{-m}}
\prod_{j=1}^{k}\prod_{i_j=1}^{w_j}(1-q^{r_j+m(i_j-1)})q^{-m(i_j-1)}\\
&\hspace{4.3cm}\times\frac{\prod_{i_j=1}^{n+u_k-u_j}(1-q^{r-s_j+m(i_j-1)})q^{r_j-m(i_j-1)}}
{\prod_{i_{j-1}=1}^{n+u_k-u_{j-1}}(1-q^{r-s_{j-1}+m(i_{j-1}-1)})q^{-m(i-1)}}.
\end{align*}
Then, proceeding as in the derivation of Theorem 3.3, the required limiting expressions (\ref{eq4.2}) and (\ref{eq4.3}) are deduced.
\end{document}